\documentclass[12pt]{article}
\usepackage{amssymb}
\usepackage{graphicx}
\usepackage{amsmath}

\RequirePackage[latin1]{inputenc} \RequirePackage[T1]{fontenc}
\newtheorem{theorem}{Theorem}

\newtheorem{lemma}{Lemma}[section]

\newenvironment{proof}[1][Proof]{\textbf{#1.} }{\ \rule{0.5em}{0.5em}}
\def \N{\mbox{I\hspace{-.15em}N}}
\def \R{\mbox{I\hspace{-.15em}R}}
\def \P{\mbox{I\hspace{-.15em}P}}
\def \E{\mbox{I\hspace{-.15em}E}}

\begin{document}

\author{Auguste Aman\thanks{E-mail address: augusteaman5@yahoo.fr (Corresponding author)}\, and Modeste N'Zi\thanks{E-mail address: modestenzi@yahoo.fr} \\
%EndAName
UFR de Math\'{e}matiques et Informatique,\\
22 BP 582 Abidjan 22,\ C\^{o}te d'Ivoire\\
}\date{}
\title{Homogenization of reflected semilinear  PDE with nonlinear
Neumann boundary condition.} \maketitle
\begin{abstract}
We study the homogenization problem of semi linear reflected partial
differential equations (reflected PDEs for short) with nonlinear
Neumann conditions. The  non-linear term is a function of the
solution but not of its gradient. The proof are fully probabilistic
and uses  weak convergence of associated reflected generalized backward differential
stochastic equations (reflected GBSDEs in short).
\end{abstract}

\section{Introduction}
Let $L^{\varepsilon }_{x}$ be an uniformly elliptic second-order partial differential operator
indexed by a parameter $\varepsilon> 0$. The homogenization
problem consists in analyzing the behavior as $\varepsilon \downarrow 0$, of the
solution $u_{\varepsilon}$ for the partial differential equations (PDEs): $L^{\varepsilon}_x
u_{\varepsilon}=f$ in a domain $\Theta$ of $\R^{d}$. The coefficients of the operator $L^{\varepsilon }_{x}$ are subject to an
appropriate conditions. For example,  periodicity by Freidlin \cite{F1}, almost periodicity by Krylov \cite{Ko}. In this paper we shall consider the case of locally periodic coefficients introduced by Ouknine et al. \cite{Oual}. In the probabilistic approach the problem becomes the following. What is the limit of the laws of the diffusion processes
$X^{\varepsilon}$ with generator $L^{\varepsilon}_{x}$, as $\varepsilon\downarrow
0$ ? This problem has been studied for diffusion processes in the whole of
$\R^{d}$ by Krylov\cite{Ko}, Freidlin\cite{F1} and Bensoussan\cite{BLP}. In the case of the presence of boundary
conditions, we can see the work of Tanaka \cite{T} and recently this one of Ouknine et al.\cite{Oual}.
In \cite{P1}, Pardoux has combined the probabilistic approach of chapter 3 of \cite{BLP} to derive
homogenization results for semi linear parabolic partial
differential equation involving a second order differential operator
of parabolic type in the whole space with highly oscillating drift
and nonlinear term. Furthermore Pardoux and Ouknine \cite{OP} pursue
the same program for a semi-linear elliptic PDE with nonlinear
Neumann boundary conditions.

This paper is devoted to prove a
similar result for obstacle problem of semi-linear elliptic partial
differential equations (PDEs in short) with a nonlinear Neumann boundary conditions $g$ satisfying moreover linear growth conditions. We
use for this purpose the link between backward stochastic differential
equations (BSDEs in short) and PDEs.

The concrete formulation of BSDEs in the nonlinear case was first introduced by Pardoux
and Peng $\cite{PP1}$ who proved existence and uniqueness of adapted
solutions for these equations, under suitable conditions on the
coefficient and the terminal value. They provide probabilistic
formulas for solutions of semi linear PDEs which generalized the well-known Feynman-Kac formula. The interest of this
kind of stochastic equations has increased steadily, since it has
been widely recognized that they provide a useful framework
for the formulation of many problems in mathematical finance (see
$\cite{EPQ})$, in stochastic
control and differential games (see $\cite{HL}$). A class of
BSDEs called generalized BSDEs has been introduced by Pardoux and Zhang in $\cite{PZ}$. This kind of BSDEs has an additional integral with respect to a continuous increasing
process and is used to provide probabilistic formulas for solutions of semi linear PDEs with nonlinear Neumann boundary conditions. Further, in order to provide probabilistic formula for solutions of obstacle problem for PDEs with nonlinear Neumann boundary conditions, Ren and Xia \cite{Ral} considered reflected generalized BSDEs.
Our goal in this paper is to give a homogenization result of obstacle problem for semi-linear PDEs with nonlinear Neumann boundary conditions via a weak convergence of reflected generalized BSDEs.

To describe our result more precisely, we first recall some
classical notations that will be used in the sequel of this paper.
Let $\mathcal{C}^{1}_{b}(\R^{d})$ be the space of real continuously differentiable functions such that the derivatives are bounded and $\Theta$ be a regular convex and bounded  domain in $\R^{d}$. We introduced a function $\rho\in\mathcal{C}^{1}_{b}(\R^{d})$ such
that $\rho=0$ in $\bar{\Theta}$,\, $\rho>0$ in $\R^{d}\backslash
\bar{\Theta}$ and $\rho(x)=(d(x,\bar{\Theta}))^{2}$ in the
neighborhood of $\bar{\Theta}$. Let us note that, since the domain
$\Theta$ is smooth, it is possible to consider an extension
$\psi\in\mathcal{C}^{2}_{b}(\R^{d})$ of the function
$d(x,\partial\Theta)$ defined on the restriction to $\Theta$ of the
neighborhood of $\partial\Theta$ such that $\Theta$ and
$\partial\Theta$ are characterized by
\begin{eqnarray*}
\Theta=\left\{x\in\R^{d}:\;\; \psi(x)>0\right\}\quad\quad\mbox{and}\qquad\qquad
\partial\Theta=\left\{x\in\R^{d}:\;\; \psi(x)=0\right\}
\end{eqnarray*}
and for all $x\in\partial\Theta$, $\nabla\psi(x)$ coincide with the unit normal pointing toward the interior of
$\Theta$ (see for example $\cite{LS}$,\;Remark $3.1$). In particular we may and do choose $\rho$ and $\psi$ such that
$\langle \nabla \psi(x),\delta(x)\rangle \leq 0,$ for all $x\in\R^{d}$, where $\delta(x)=\nabla \rho(x)$ which we call the penalization term. Now let us consider the following reflected semi-linear
partial differential equations with Neumann boundary condition for $\varepsilon>0)$:\\\\
$
\left\{
\begin{array}{l}
\min\left(u^{\varepsilon}(t,x)-h(t,x),\right.\\\\
\left.-\frac{\partial u^{\varepsilon}}{\partial
t}(t,x)-L^{\varepsilon}_{x}u^{\varepsilon}-f(x,u^{\varepsilon}(t,x))
\right)=0,\; if\, (t,x)\in\R_{+}\times\Theta,\\\\
\Gamma^{\varepsilon}_{x}u^{\varepsilon}(t,x)+g(x,u^{\varepsilon}(t,x))=0,\;\; if\, (t,x)\in \R_{+}\times\partial \Theta,\\\\
u^{\varepsilon}(0,x)=l(x),\;\;\;  x\in\overline{\Theta},
\end{array}\right.
$
\\\\
where $\Gamma^{\varepsilon}_{x}$ is defined in $(\ref{a2})$.
Under suitable linear growth, monotonic and Lipschitz conditions on
the coefficients $f,\; g$ and $l$, we will show that
$u^{\varepsilon}(t,x)$ the viscosity solution of the previous reflected PDEs
converge, as $\varepsilon$ goes to $0$, to a function $u(t,x)$, the  viscosity solution to the following reflected PDEs
with nonlinear Neumann boundary conditions:\\\\
$
\left\{
\begin{array}{l}
\min\left(u(t,x)-h(t,x),\right.\\\\
\left.-\frac{\partial u}{\partial t}(t,x)-L^{0}_{x}u-f(x,u(t,x))\right)=0,\;\, if \; (t,x)\in\R_{+}\times\Theta,\\\\
\Gamma^{0}_x u(t,x)+g(x,u(t,x))=0,\;\; \; if\; (t,x)\in \R_{+}\times\partial \Theta,\\\\
u(0,x)=l(x),\;\;\; x\in\overline{\Theta},
\end{array}\right.
$
\\\\
where operators $L^{0}_x$ and $\Gamma^{0}_x$ are in the form $(\ref{op})$.

The rest of this paper is organized as follows. In section $2$ we recall
without proof, some results due to Ouknine et al.\cite{Oual} on the
homogenization of the reflected diffusion in the bounded domain with
locally periodic coefficients and some basic notations on generalized BSDEs. Section $3$ is devoted to the statement and proof of our main result. Finally, in Section 4, we give a homogenization result for some PDEs with nonlinear Neumann boundary conditions.

\section {Locally periodic homogenization of reflected diffusion}
\setcounter{theorem}{0} \setcounter{equation}{0}
In this section, we state some results which were proved in Ouknine et al. \cite{Oual}, and which we will use in the rest of the paper. Let us recall that $\Theta$ is a regular convex and bounded  domain in $\R^{d}; (d\geq
1)$ defined as previous and let
\begin{eqnarray}
L^{\varepsilon }_{x} &=&\sum_{i,j}^{d}a_{ij}(x,\frac{x}{\varepsilon
})\partial _{i}\partial _{j}+\frac{1}{\varepsilon }\sum_{i=1}^{d}
b_{i}(x,\frac{x}{\varepsilon } )\partial _{i}+\sum_{i=1}^{d}
c_{i}(x,\frac{x}{\varepsilon }
)\partial _{i}\label{a2.}
\end{eqnarray}
and
\begin{eqnarray}
\Gamma^{\varepsilon
}_{x}&=&\sum_{i=1}^{d}\partial_{i}\psi(x,\frac{x}{\varepsilon
})\partial _{i} \label{a2}
\end{eqnarray}
be given (with $\varepsilon>0$), where $\partial_{i}=\partial/\partial x_i$. We also set
\begin{eqnarray*}
L_{x,y} &=&\sum_{i,j}^{d}a_{ij}(x,y)\partial _{i}\partial _{j}+\sum_{i=1}^{d}
b_{i}(x,y)\partial _{i},\;\; x,y\in\R^{d}
\end{eqnarray*}
and require the following:
\begin{description}
\item ({\bf H.1}) $L_{x,y}$ is uniformly elliptic and the matrix $a(x)=[a_{i,j}]$ is factored as $\sigma(x,y)\sigma^{*}(x,y)/2$
\item ({\bf H.2}) The functions $\sigma:
\R^{d}\times\Theta\rightarrow \R^{d},\; b:
\R^{d}\times\Theta\rightarrow \R^{d},\;c:
\R^{d}\times\Theta\rightarrow \R^{d}$ are locally periodic (i.e,
periodic with respect to second variable; of period $1$ in each
direction in $\Theta$).\\\\
$\left\{
\begin{array}{l}
(i)\quad \mbox{{\it Global Lipschitz condition}: there exists a constant}\; C \;\mbox{such that for any}\, \zeta = \sigma,b, c,
\\
\|\zeta(x, y)-\zeta(x', y')\|\leq C\left(\|x-x'\|+\|y- y'\|\right),\; \forall\;x,x'\in\R^{d},\;  y, y' \in\Theta,
\\\\
(ii) \quad\mbox{The partial derivatives} \;\partial_{x}\zeta(x, y)\;\mbox{as well as the mixed derivatives}\;
\partial_{xy}^{2}\zeta(x, y)\;\\ \mbox{exist and are continuous},\; \forall x \in\R^{d}, y \in\Theta\;\mbox{for}\;\zeta = a,b, c
\\\\
(iii)\quad\mbox{The coefficients are bounded, that is, there exists a constant}\; C \;\mbox{such that for any}\\
\zeta = a,b, c\qquad
\|\zeta(x, y)\|\leq C,\;\; x \in \R^{d}, y\in\Theta.
\end{array}\right.
$
\end{description}

The differential operator $L^{\varepsilon}_{x}$ inside $\Theta$ together with the boundary condition $\Gamma_{x}^{\varepsilon}u=0$ on $\partial\Theta$ determine a unique diffusion process $X^{\varepsilon}$ in $\bar{\Theta}$ which we call the $(L^{\varepsilon}_{x},\Gamma^{\varepsilon}_x)$-diffusion. Given a $d$-dimensional Brownian motion $(W_{t}: \; t\geq 0)$
defined on a complete probability space $(\Omega,\mathcal{F},\P)$, let $(X^{\varepsilon}, G^{\varepsilon})$ be the unique solution with values in $\Theta\times\R_{+}$ of the following reflected SDE:
\begin{eqnarray}
\left\{
\begin{array}{l}
dX^{\varepsilon}_{t}=\frac{1}
{\varepsilon}b(X^{\varepsilon}_{t},\frac{X^{\varepsilon}_{t}}{\varepsilon})dt
+c(X^{\varepsilon}_{t},\frac{X^{\varepsilon}_{t}}{\varepsilon})dt
+\sigma(X^{\varepsilon}_{t},\frac{X^{\varepsilon}_{t}}{\varepsilon})dW_{t}
+\nabla\psi(X^{\varepsilon}_{t})dG_{t}^{\varepsilon},\ t \geq 0
\\\\
X^{\varepsilon} _{t}\in \Theta,\ G^{\varepsilon}\; \mbox{is continuous
and increasing such that},\\\\
\int_{0}^{t}\nabla\psi(X^{\varepsilon}_{s}){\bf 1}_{\{X^{\varepsilon}_{s}\in\Theta\}}dG_{s}^{\varepsilon}=0
,\ t \geq 0  \\\\
X^{\varepsilon} _{0}=x.
\end{array}
\right.    \label{a4}
\end{eqnarray}
By requirement there exists a $L_{x,y}$-diffusion on $\R^{d}$ with generator $L_{x,y}$ and by $y$-periodicity
assumption on the coefficients this process induces diffusion process $U^{x}$ on the
$d$-dimensional torus ${\bf T}^{d}$, moreover this diffusion process is ergodic. We denote by $m(x,.)$
its unique invariant measure. In order for the process with generator $L^{\varepsilon}_{x}$
to have a limit in law as $\varepsilon\downarrow 0$, we need the following condition to be in force.
\begin{description}
\item $({\bf H.3})$ {\it Centering condition}: for all $x$,
\begin{eqnarray*}
\int_{\Theta}b(x,u)m(x,du) = 0.
\end{eqnarray*}

\end{description}
Let us denote
\begin{eqnarray*}
\widehat{b}=
\left(
\begin{array}{c}
\widehat{b}_{1} \\
\widehat{b}_{2 }\\
. \\
.  \\
. \\
. \\
\widehat{b}_{d}\\
\end{array}
\right)
\end{eqnarray*}
with $\displaystyle{\widehat{b}^{k}(x,u) =\int^{\infty}_{0}\E_{u}\{b^{k}(x,U^{x}_{
t} )\}dt}$, where under $\P_{u},\, U^x$ starts from $u$.

We put
\begin{eqnarray}
L^{0}_x=\sum_{i,j=1}^{d}\overline{A}_{0}^{i j}(x)\partial_{x_{i}x_{j}} +\sum_{i=1}^{d}\overline{C}_{0}^{i}(x)\partial_{x_{i}},\;\; \Gamma^{0}_x=\sum_{i=1}^{d}\gamma^{0}_{i}(x)\partial_{x_{i}},\label{op}
\end{eqnarray}
where the coefficients are respectively defined by
\begin{align*}
C_{0}(x, y) = \left[\partial_{x}\widehat{b}b+(I+\partial_{y}\widehat{b})c+\frac{1}
{2}Tr\partial^{2}_{xy}\widehat{b}\sigma\sigma^{*}\right](x,y),\\\\
S_{0}(x, y) = \left[(I+\partial_{y}\widehat{b})\sigma\right],\;\; A_{0}(x, y) = S_{0}S_{0}^{*} (x, y),\\
\overline{C}_{0}(x) = \int_{{\bf T}^{d}}C_{0}(x,u)m(x,du),\; \overline{A}_{0}(x)=\int_{{\bf T}^{d}}A_{0}(x,u)m(x,du),\\\\
\gamma^{0}(x)=\int_{{\bf T}^{d-1}}\left[(I+\partial_{y}\widehat{b})\nabla\psi\right](x, y)d\widetilde{m}(x,dy).
\end{align*}

Let us consider the following SDE:\\\\
$
\left\{
\begin{array}{l}
dX_{t}=\overline{A}^{1/2}_{0}(X_{t})dW_{t}+\overline{C}_{0}(X_{t})dt+\gamma^{0}(X_{t})dG_{t},\; t\geq 0,\\\\
X\in\Theta,\; G \,\mbox{is continuous and increasing },\\\\
\int_{0}^{t}\gamma^{0}(X_{s}){\bf 1}_{\{X_{s}\in\Theta\}}dG_{s}=0,\, t\geq 0,\\\\
X_{0}=0.
\end{array}
\right.
$\\\\
The operators $L^{0}$ and $\Gamma^{0}$ are acting on $x$.

We can now state the main result of Ouknine et al. \cite{Oual} (see Theorem 3.1), which is a generalization of a result of Tanaka \cite{T} (see Theorem 2.2).
\begin{theorem}
\label{theorem1} Under assumptions $({\bf H1})$-$({\bf H3})$, the
$(L^{\varepsilon}_x,\Gamma^{\varepsilon}_x)$-diffusion process
$X^{\varepsilon}$ converges in law to an
$(L^{0}_x,\Gamma^{0}_x)$-diffusion $X$ as $\varepsilon\rightarrow 0$.\\
Moreover,
\begin{eqnarray*}
(X^{\varepsilon},M^{X^{\varepsilon}},G^{\varepsilon})\Longrightarrow (X,M^{X},G),
\end{eqnarray*}
where $M^{X}$ (resp., $M^{X^{\varepsilon}}$) is the martingale part of $X$ (resp.\,$X^{\varepsilon}$),  $G$ (resp., $G^{\varepsilon}$) is the local time of $\psi(X)$ (resp., $\psi(X^{\varepsilon})$ ) at $0$ and\\
"$\Longrightarrow$" designates the weak convergence in the space $\mathcal{C}([0,T],\R^{d})\times\mathcal{D}([0,T],\R^{d})\times\mathcal{C}([0,T],\R^{d})$ endowed with the product topology of the uniform norm for the space $\mathcal{C}([0,T],\R^{d})$ and the $S$-topology for $\mathcal{D}([0,T],\R^{d})$.
\end{theorem}

It  follows moreover easily from the result of Tanaka \cite{T} that
\begin{lemma}
Under assumptions $({\bf H1})-({\bf H3}),$
\begin{eqnarray*}
\sup_{\varepsilon}\E\left(\sup_{0\leq s\leq t}|X^{\varepsilon}_{s}|^{p}+G^{\varepsilon}_{t}\right)<\infty,
\end{eqnarray*}
for any $p\geq 1$.
\end{lemma}
\section{Reflected generalized BSDEs and weak convergence}
\setcounter{theorem}{0} \setcounter{equation}{0}
\subsection{Preliminaries and statement of the problem}
Let us introduce some spaces:
\begin{eqnarray*}
H^{2,d}(0,t)&=&\left\{\{\psi_s;\,0\leq s\leq t\}:\, \mbox{$\R^{d}$-valued predictable process}\, \right.\\
&&\left.\mbox{such that}\, \E\int_{0}^{t}|\psi(s)|^{2}ds<\infty\right\},\\\\
S^{2}(0,t)&=&\left\{\{\psi_s;\,0\leq s\leq t\}:\, \mbox{real-valued, progressively measurable process}\right.\\
&&\left.\mbox{such that}\, \E\left[\sup_{0\leq s\leq t}|\psi(s)|^{2}\right]<\infty\right\},\\\\
A^{2}(0,t)&=&\left\{\{K_s;\,0\leq s\leq t\} \mbox{real-valued adapted, continuous, },\right.\\
&&\left.\mbox{and increasing process such that}\, K(0)=0,\E|K(t)|^{2}<\infty\right\}.
\end{eqnarray*}
In addition, we give the following assumptions:

Let $f,g:\overline{\Theta}\times \R \rightarrow \R$ be a continuous
function satisfying the following assumptions:
\\ There exist constants $C>0,\, p
\geq 1,\ \mu\in \R $\, and$\ \beta<0 $, such that for  all $x\in
\bar{\Theta},\ y,\ y'\in \R,$\\\\
$(\mbox{\bf H.4})
\left\{
\begin{array}{l}
(f.i)\quad \left(y-y^{\prime }\right)\left(f(x,y)-f(x,y^{\prime
})\right) \leq \mu
|y-y^{\prime }|^{2},\\\\
(f.ii)\quad  |f(x,y)|\leq C(|x|^{p}+|y|^{2}),\\\\
(g.i) \quad(y-y^{\prime})\left(g(x,y)-g(x,y^{\prime })
\right) \leq \beta |y-y^{\prime }|^{2}, \\\\
(g.ii) \quad |g(x,y)|\leq C(|x|^{p}+|y|^{2}).
\end{array}\right.
$
\\
Moreover, let $l:\Theta \rightarrow \R$ be a continuous function and $h\in C^{1,2}(\R_{+}\times \overline{\Theta})$ for which there exist a constant $C$ such that\\\\
${\bf H.5} \left\{
\begin{array}{l}
(i)\ l(x)\leq C(1+|x|^{p}),\\\\
(ii)\ h(t,x) \leq C(1+| x|^{p}),\\\\\
(iii)\ h(0,x)\leq l(x).
\end{array}\right.
$\\
where $p$ is obtained in $({\bf H4})$.

For each fixed $(t,x)\in\R_{+}\times\overline{\Theta}$, we consider
$(Y_{s}^{\varepsilon },Z_{s}^{\varepsilon },K_{s}^{\varepsilon })_{0\leq s\leq t}$, the unique solution of
reflected generalized BSDE
\begin{eqnarray}
\left\{
\begin{array}{l}
Y_{s}^{\varepsilon }=l(X_{t}^{\varepsilon})+\int_{s}^{t}
f(X_{r}^{\varepsilon},Y_{r}^{\varepsilon })dr
+\int_{s}^{t}g(X_{r}^{\varepsilon},Y_{r}^{\varepsilon
})dG_{r}^{\varepsilon}+K_{t}^{\varepsilon}-K_{s}^{\varepsilon }
-\int_{s}^{t}Z_{r}^{\varepsilon}dM^{X^{\varepsilon}}_{r},
\\\\
Y_{s}^{\varepsilon}\geq h(s,X_{s}^{\varepsilon }),\\\\
K^{\varepsilon }\; \mbox{is nondecreasing such that} \
K_{0}^{\varepsilon }=0\  \mbox{and}\ \int_{0}^{t}(Y_{s}^{\varepsilon
}-h(s,X_{s}^{\varepsilon }))dK_{s}^{\varepsilon }=0.
\end{array}
\right. \label{a6}
\end{eqnarray}
Let us note that existence and uniqueness of a solution to $(\ref{a6})$ is proved in \cite{Ral}.

Furthermore, we also consider $(Y_{s},Z_{s},K_{s})_{0\leq
s\leq t}$ the unique solution of the reflected BSDE
\begin{eqnarray}
\left\{
\begin{array}{l}
Y_{s}=l(X_{t})+\int_{s}^{t} f(X_{r},Y_{r})dr +\int_{s}^{t}g(X_{r},Y_{r})dG_{r}-\int_{s}^{t}Z_{r}dM^{X}_{r}+K_{t}-K_{s} \\\\
Y_{s}\geq h(s,X_{s})\\\\
K \;\mbox{is nondecreasing such that} \, K_{0}=0\; \mbox{and}\;
\int_{0}^{t}(Y_{s}-h(s,X_{s}))dK_{s}=0.
\end{array}\right.
\label{Eq1}
\end{eqnarray}
Throughout the rest of this paper, we put: \[ M_{s}^{\varepsilon
}=\int_{0}^{s}Z_{r}^{\varepsilon}dM^{X^{\varepsilon}},\;
M_{s}=\int_{0}^{s}Z_{r}dM^{X}_{r};\;  0\leq s\leq t\] and we consider
$(Y,M,K)$ (resp.,\
$(Y^{\varepsilon},M^{\varepsilon},K^{\varepsilon})$)
as a random element of the space $D([0,t];\R)\times
D([0,t];\R)\times C([0,t];\R_{+})$. $D([0,t];\R)$ denotes the Skorohod space (the space of "cadlag" functions), endowed with the so-called $S$-topology
of Jakubowski \cite{J} and $C([0,t];\R_{+})$ the space of functions of $[0, t]$ with values in $\R^d$ equipped with the topology of
uniform convergence.
\subsection{Weak convergence}
The goal of this section is to prove the following result
\begin{theorem}
\label{C5T1} Under conditions $(\bf H1)$-$(\bf H4)$ the family  $(Y^{\varepsilon },M^{\varepsilon },K^{\varepsilon })_{\varepsilon>0} $ of
processes converge in law to $(Y,M,K) $ on the space $D([ 0,t],\R)\times D([ 0,t],\R)\times
C([ 0,t];\R_{+}).$ as $\varepsilon\rightarrow 0$.\\
Moreover, $$\lim_{\varepsilon\rightarrow
0}Y^{\varepsilon}_{0}=Y_{0}\;\; \mbox{in} \;\; \R.$$
\end{theorem}
The proof of this theorem relies on four
lemmas. The first one is a simplified version of Theorem 4.2 in  (Billingsley \cite{BI},
p. $25)$ so we omit the proof. We note also that the convergence is taken in law.
\begin{lemma}
\label{C5P1}
Let $\left( U^{\varepsilon }\right) _{\varepsilon>0}$\ be a family
of random variables defined on the same probability space. For
each  $\varepsilon \geq 0$,\ we assume the existence of a family
random variables $\left( U^{\varepsilon ,n}\right) _{n\leq 1}$,\ such
that:
\begin{description}
\item $(i)$ For each \,$n\geq 1,\;\ U^{\varepsilon ,n}\Longrightarrow U^{0,n}$ as
$\varepsilon$ goes to z\'{e}ro and
\item $(ii)$ $U^{0,n}\Longrightarrow U^{0}$ as $n\longrightarrow +\infty$.
\item $(iii)$ We suppose further that\, $U^{\varepsilon
,n}\ \Longrightarrow U^{\varepsilon }$ as $
n\longrightarrow +\infty$,\  uniformly in $\varepsilon $
\end{description}
Then, $U^{\varepsilon }$ converges in law to $U^{0}$.
\end{lemma}
For the last third lemmas, let us introduce the penalization method to construct a sequence of GBSDE.
For each $n\in \mathbb{N}^{\ast},$ we set
\begin{equation*}
f_{n}(x,y)=f(x,y)+n(y-h(s,x))^{-}
\end{equation*}
Thanks to the result of Pardoux and Zhang \cite{PZ}, for each $n\in
\N^{\ast }$, there exists a unique pair of
$\mathcal{F}_{t}-$progressively measurable processes
$\left(Y^{\varepsilon,n},Z^{\varepsilon ,n}\right)$  with values in
$ \R\times \R^{d}$ satisfying
\begin{eqnarray*}
\E\left[\sup_{0\leq s\leq t}|Y^{\varepsilon,n}|^{2}+\int_{0}^{t}|Z^{\varepsilon,n}_s|^{2}ds\right]<\infty
\end{eqnarray*}
and
\begin{eqnarray}
Y_{s}^{\varepsilon ,n} &=&l(X_{t}^{\varepsilon
})+\int_{s}^{t}f_{n}(X_{r}^{\varepsilon },Y_{r}^{\varepsilon
,n})dr+\int_{s}^{t}g(X_{r}^{\varepsilon },Y_{r}^{\varepsilon
,n})dG_{r}^{\varepsilon }  \nonumber \\
&&-(M_{t}^{\varepsilon ,n}-M_{s}^{\varepsilon ,n}),  \label{T1}
\end{eqnarray}
where
\[
M^{\varepsilon ,n}_{s}=\int_{0}^{s}Z_{r}^{\varepsilon ,n}dM^{X^{\varepsilon}}_{r}.
\]

From now on, $C$ is a generic constant that may vary from line to another.
\begin{lemma}
\label{C5L1}Assume $({\bf H1})$-$({\bf H4})$. Then for any $n$, there exists the unique process $ (Y^{0,n},M^{0,n})$ such that the family of processes$\left( Y^{\varepsilon ,n},M^{\varepsilon ,n}\right)_{\varepsilon>0}$ converges to it in $D([0,T],\R)\times D([0,T],\R)$, as $\varepsilon$ goes to $0$.
\end{lemma}
\begin{proof}
{\bf Step 1. A priori estimate  uniformly in $\varepsilon $} and $n$
\smallskip
\newline Applying It\^{o}'s formula to the function
$|Y_{t}^{\varepsilon ,n}|^{2}$, we get
\begin{eqnarray}
|Y_{s}^{\varepsilon ,n}|^{2}+\int_{s}^{t}|Z_{r}^{\varepsilon
,n}|^{2}dr &=&|l(X_{t}^{\varepsilon })|^{2}+2\int_{s}^{t}
Y_{r}^{\varepsilon ,n}f\left(X^{\varepsilon }_{r},Y_{r}^{\varepsilon
,n}\right)  dr
 \nonumber\\
&&+2\int_{s}^{t} Y_{r}^{\varepsilon ,n}g\left(X_{r}^{\varepsilon
},Y_{r}^{\varepsilon ,n}\right)  dG_{r}^{\varepsilon }\nonumber \\
&&-2\int_{s}^{t} Y_{r}^{\varepsilon ,n}Z_{r}^{\varepsilon ,n}dM^{X^{\varepsilon}}_{r}
+2\int_{s}^{t} Y_{r}^{\varepsilon ,n}dK_{r}^{\varepsilon ,n},
\label{5.12}
\end{eqnarray}
where
\begin{eqnarray}
K_{s}^{\varepsilon ,n}=n\int_{0}^{s}[Y_{r}^{\varepsilon
,n}-h(r,X^{\varepsilon ,n}_{r})]^{-}dr
\label{Ra}
\end{eqnarray}
 By using assumptions
$({\bf H3})$ and $({\bf H4}) $, the inequality $ \int_{s}^{t}
Y_{r}^{\varepsilon ,n}dK_{r}^{\varepsilon ,n} \leq \int_{s}^{t}
h\left(r,X^{\varepsilon }_{r}\right)dK_{r}^{\varepsilon ,n} $ and
Young's inequality, we get for every $\gamma>0$
\begin{eqnarray*}
\E|Y_{s}^{\varepsilon ,n}|^{2}+\int_{s}^{t}|Z_{r}^{\varepsilon
,n}|^{2}dr &\leq &\left(1+C\E|X_{t}^{\varepsilon }|^{2p}\right)+\E%
\int_{s}^{t}\left\{ \left( 2\mu +1\right) |Y_{r}^{\varepsilon
,n}|^{2}+C^{2}|X_{r}^{\varepsilon }|^{2q}\right\} dr \\
&&+\E\int_{s}^{t}\left\{(2\beta +\gamma)\left| Y_{r}^{\varepsilon
,n}\right| ^{2}+\frac{2C^{2}}{\gamma}\left|
X_{r}^{\varepsilon }\right| ^{2q}\right\}dG_{r}^{\varepsilon } \\
&&+\E\int_{s}^{t} h\left(r,X^{\varepsilon
}_{r}\right)dK_{r}^{\varepsilon ,n} .
\end{eqnarray*}
By virtue of assumption $(h.i)$, Young's equality and choosing $\gamma=|\beta|$ we obtain for every $\delta>0$
\begin{eqnarray*}
\E\left(\left| Y_{s}^{\varepsilon ,n}\right| ^{2}+\int_{s}^{t}\left|
Z_{r}^{\varepsilon ,n}\right| ^{2}dr+\left| \beta \right|
\int_{s}^{t}\left| Y_{r}^{\varepsilon ,n}\right|
^{2}dG_{r}^{\varepsilon }\right) &\leq &C(1+\E\int_{s}^{t}\left|
Y_{r}^{\varepsilon ,n}\right| ^{2}dr)
\nonumber \\
&&+\frac{1}{\delta }\E\left(\sup_{0\leq s\leq t}\left|
h\left(s,X^{\varepsilon
}_{s}\right)\right| ^{2}\right)  \nonumber \\
&&+\delta \E\left(K_{t}^{\varepsilon ,n}-K_{s}^{\varepsilon
,n}\right)^{2}.
\end{eqnarray*}
But
\begin{eqnarray*}
\E(K_{t}^{\varepsilon ,n}-K_{s}^{\varepsilon ,n})^{2} &\leq
&C\E\left(|Y_{s}^{\varepsilon ,n}|^{2}+|l(X_{t}^{\varepsilon
})|^{2}+\int_{s}^{t}(|X_{r}^{\varepsilon }|^{2q}+|Y_{r}^{\varepsilon
,n}|^{2})dr\right)   \nonumber \\
&&+C\E\int_{s}^{t}(|X_{r}^{\varepsilon }|^{2q}+|Y_{r}^{\varepsilon
,n}|^{2})dG_{r}^{\varepsilon }.
\end{eqnarray*}
So, if $\delta =\inf \left\{ \dfrac{1}{2C},\dfrac{|\beta
|}{2C}\right\},$
\[
\E|Y_{s}^{\varepsilon ,n}|^{2}+\E\int_{s}^{t}|Z_{r}^{\varepsilon
,n}|^{2}dr+\E\int_{s}^{t}|Y_{r}^{\varepsilon
,n}|^{2}dG_{r}^{\varepsilon }\leq
C\left(1+\E\int_{s}^{t}|Y_{r}^{\varepsilon ,n}|^{2}dr\right)
\]
and it follows from the Gronwall's lemma that
\begin{eqnarray*}
\E\left(|Y_{s}^{\varepsilon ,n}|^{2}+\int_{s}^{t}|Z_{r}^{\varepsilon
,n}|^{2}dr+\int_{s}^{t}|Y_{r}^{\varepsilon
,n}|^{2}dG_{r}^{\varepsilon }+|K_{t}^{\varepsilon
,n}|^{2}\right)&\leq C,\;\;\;\;0\leq s\leq t,\;\; \forall\ n\in
\N^{*}, \forall\ \varepsilon>0.
\end{eqnarray*}
Finally, Burkhölder-Davis-Gundy inequality applying to (\ref{5.12})
gives the following:
\begin{eqnarray*}
\sup_{n\in\N^{*}}\sup_{\varepsilon>0}\E\left(\sup_{0\leq s\leq
t}|Y_{s}^{\varepsilon ,n}|^{2}+\int_{0}^{t}|Z_{r}^{\varepsilon
,n}|^{2}dr+\int_{0}^{t}|Y_{r}^{\varepsilon
,n}|^{2}dG_{r}^{\varepsilon }+|K_{t}^{\varepsilon
,n}|^{2}\right)\leq C.
\end{eqnarray*}

{\bf Step 2: Tightness }

We adopt the point of view of the $S$-topology of Jakubovski
\cite{J}. In fact we define the conditional variation of the process
$Y^{\varepsilon ,n}$ on the interval $[0,t]$ as the quantity
\begin{eqnarray*}
CV_{t}(Y^{\varepsilon
,n})=\sup_{\pi}\E\left(\sum_{i}|\E(Y_{t_{i+1}}^{\varepsilon
,n}-Y_{t_{i}}^{\varepsilon
,n}/\mathcal{F}_{t_{i}}^{\varepsilon})|\right)
\end{eqnarray*}
where the supremum is taken over all partitions $\pi:0=t_0<t_1<...<t_m=t$ of the interval
 $[0,t]$.

 Clearly
\begin{eqnarray*}
CV_{t}(Y^{\varepsilon ,n})\leq \E\left(\int_{0}^{t}|f(X_{r}^{%
\varepsilon },Y_{r}^{\varepsilon
,n})|dr+\int_{0}^{t}|g(X_{r}^{\varepsilon },Y_{r}^{\varepsilon
,n})|dG_{r}^{\varepsilon }+|K_{t}^{\varepsilon ,n}|\right),
\end{eqnarray*}
and it follows from step $1$ and  $({\bf H3})$ that
\begin{eqnarray*}
\sup_{\varepsilon }\left(CV_{t}(Y^{\varepsilon
,n})+\E\left(\sup_{0\leq s\leq t}|Y_{s}^{\varepsilon
,n}|+\sup_{0\leq s\leq t}|M^{\varepsilon
,n}_{s}|\right)+|K^{\varepsilon,n}_{t}|^{2}\right)<+\infty.
\end{eqnarray*}
Let us denote
\begin{eqnarray*}
H_{s}^{\varepsilon
,n}=\int_{0}^{s}\left(g(X^{\varepsilon}_{r},Y^{\varepsilon
,n}_{r})+n[Y^{\varepsilon
,n}_{r}-h(r,X^{\varepsilon}_{r})]^{-}\right)dG^{\varepsilon }_{r}.
\end{eqnarray*}
By virtue of step 1 and  $({\bf H3})$, we have
\begin{eqnarray*}
\sup_{\varepsilon }\left(CV_{t}(H^{\varepsilon
,n})+\E\left(\sup_{0\leq s\leq t}|H_{s}^{\varepsilon
,n}|\right)\right)<+\infty.
\end{eqnarray*}
Hence the sequence $\left\{(Y^{\varepsilon ,n}_{s},H^{\varepsilon
,n}_{s},M^{\varepsilon ,n}_{s});\ 0\leq s \leq t \right\}$ satisfy
Jakubowski tightness criterion for quasi-martingales under
${\P}$(see \cite{J}), which is the same criterion as the one of
Meyer-Zheng (see \cite{MZ}).

{\bf Step 3: Passage to the limit and it identification}.

After extracting a subsequence which we
omit by abuse of notation, it follows by Jakubowski
tightness criterion that $\left(Y^{\varepsilon
,n},H^{\varepsilon ,n},M^{\varepsilon ,n}\right)$ converges in the weak sense to $\left(
Y^{n},H^{n},M^{n}\right)$ as $\varepsilon$ goes to 0 in $(D([ 0,t] ,\R))^{3}$ endowed with the product of
$S$-topology.\newline
Since $f$ is continuous, the mapping
$(x,y)\mapsto\int_{s}^{t}\left(f(x(r),y(r))+n[y(r)-h(r,x(r))]^{-}\right)dr$
is continuous from $\mathcal{C}([0,t],\R^{d})\times D([ 0,t] ,\R)$,
equipped respectively by topology of uniform norm and the
$S-$topology. Furthermore, it follows (see Boufoussi et al.\cite{Bal}) that for all $0\leq s\leq t$
\begin{eqnarray*}
H_{s}^{n}=\int_{0}^{s}g(X_{r},Y^{0,n}_{r})dG_{r}
\end{eqnarray*}
and
\begin{eqnarray*}
Y_{s}^{0,n} &=&l\left( X_{t}\right)
+\int_{s}^{t}f(X_{r},Y_{r}^{0,n})dr+\int_{s}^{t}g(X_{r},Y_{r}^{0,n})
dG_{r} \\
&&+n\int_{s}^{t}(Y_{r}^{0,n}-h\left(r,X_{r}\right))^{-}dr +M_{t}^{0,n}-M_{s}^{0,n}.
\end{eqnarray*}
By using a similar argument as in Pardoux \cite{P1}, one can prove that $M^{0,n}$ and $M^{X}$ are
$\mathcal{F}_{s}^{X,Y^{n}}$-martingales. Let $\left\{
(\overline{Y}_{s}^{n},\overline{U}_{s}^{n}),0\leq s\leq t\right\} $
denote the unique solution of the generalized BSDE
\begin{eqnarray*}
\overline{Y}_{s}^{n} &=&l\left( X_{t}\right) +\int_{s}^{t}f(X_{r},
\overline{Y}^{n}_{r})dr+\int_{s}^{t}g(X_{r},\overline{Y}_
{r}^{n})dG_{r} \\
&&+n\int_{s}^{t}(\overline{Y}_{r}^{n}-h\left(r,X_{r}\right))^{-}dr +
\widetilde{M}_{t}^{n}-\widetilde{M}_{s}^{n},
\end{eqnarray*}
which satisfies $\displaystyle{\E\int_{s}^{t}\overline{U}_{r}^{n}d\langle M_{r}^{X}\rangle
(\overline{U} _{r}^{n})^{*}<+\infty}$. Set also $\displaystyle{\widetilde{M}_{s}^{n}=\int_{0}^{s}\overline{U}_{r}^{n}dM_{r}^{X}}$.
Since $\overline{Y}^{n}$, $ \overline{U}^{n}$ are
$\mathcal{F}_{s}^{X}
$-adapted, and $M^{X}_{s}$-the martingale part of X is a $\mathcal{F}_{s}^{X,Y^{n}}$-martingale, so is $\widetilde{M}^{n}$. It
follows from It\^{o}'s  formula for possibly discontinuous semi-
martingales that
\begin{eqnarray*}
&&\E\left| \overline{Y}_{s}^{n}-Y_{s}^{n}\right| ^{2}+\E
[M^{n}-\widetilde{M}^{n}]_{t}-E[M^{n}-\widetilde{M}^{n}]_{s} \\
&=&2\E\int_{s}^{t}(\overline{Y}
_{r}^{n}-Y_{r}^{n})(f(X_{r},Y_{r}^{n})-f(X_{r},\overline{Y}^{n}_{r}))
dr
\\
&&+2\E\int_{s}^{t}( \overline{Y}%
_{r}^{n}-Y_{r}^{n})(g(X_{r},Y_{r}^{n})-g(X_{r},
\overline{Y}_{r}^{n}))dG_{r} \\
&&+2\E\int_{s}^{t}(\overline{Y}_{r}^{n}-Y_{r}^{n})(d\overline{K}%
^{n}-dK_{r}^{n})  \\
&\leq &\left( 2\mu +1\right) \E\int_{s}^{t}\left| \overline{Y}%
_{r}^{n}-Y_{r}^{n}\right| ^{2}dr+2\beta \E\int_{s}^{t}\left|
 \overline{Y}_{r}^{n}-Y_{r}^{n}\right| ^{2}dG_{r} \\
&\leq &C\E\int_{s}^{t}\left| \overline{Y}_{r}^{n}-Y_{r}^{n}\right|
^{2}dr).
\end{eqnarray*}
(We use the fact that $\beta<0$).\newline
Hence, from Gronwall lemma  $\overline{Y}_{s}^{n}=Y_{s}^{n}$ and $%
M_{s}^{n}=\widetilde{M}_{s}^{n},$ $0\leq s\leq t.$
\end{proof}
\begin{lemma}
\label{C5L3}Assume $\left({\bf H1}\right)$-$\left( {\bf H4}\right)$.
Then the sequence of processes $ \left( Y^{0,n},M^{0,n},K^{0,n}\right)
_{n\geq 1}$ converges in law to $ \left( Y,M,K\right) $ in
$D([0,t] ,\R)\times D([0,t] ,\R)\times C([0,t],\R_{+})$ endowed with the previous topology.
\end{lemma}

\begin{proof}
By the same computation  of step 1 and step 2 in the proof of Lemma $\ref{C5L1},$ we get
\begin{equation*}
\sup_{n}\E\left(\sup_{0\leq s\leq t}\left|
Y_{s}^{0,n}\right| ^{2}+\int_{0}^{t}\left|
Z_{r}^{0,n}\right| ^{2}dr+\int_{0}^{t}\left|
Y_{r}^{0,n}\right| ^{2}dG_{r}+\left|
K_{t}^{0,n}\right| ^{2}\right)<+\infty \label{.13}
\end{equation*}
and
\begin{eqnarray*}
\sup_{\varepsilon}\left(CV_{t}(Y^{0,n})+CV_{t}(K^{0,n})+\E\left(\sup_{0\leq s\leq t}|M^{0,n}_{s}|^{2}\right)\right)< +\infty.
\end{eqnarray*}
Further $\left(Y^{n,0},K^{0,n},M^{0,n}\right)$ satisfy
Jakubowski tightness criterion for quasi-martingales under
${\P}$ and by the similar argument as in Lemma 3.2, $\left(Y^{.,n},F^{0,n},M^{0,n}\right)$, converges in law to $\left(
\widehat{Y},\widehat{M},\widehat{K}\right)$ as $n$ goes to $\infty$
on $D([ 0,t] ,\R^{2})\times C([0,t],\R_{+})$ endowed with the previous topology.
Finally by the uniqueness of the solution of the reflected generalized BSDE $(\ref{Eq1})$, it follows that $\left(
\widehat{Y},\widehat{M},\widehat{K}\right)=\left(
Y,M,K\right)$ which ends the proof.
\end{proof}
\begin{lemma}
\label{C5L2}Assume $({\bf H1})$-$({\bf H4})$, the  family of processes $\left(
Y^{\varepsilon ,n},M^{\varepsilon ,n},K^{\varepsilon ,n}\right)
_{n\geq 1}$ converge uniformly in $\varepsilon \in (0,1]$ in probability
to the family of processes $\left( Y^{\varepsilon },M^{\varepsilon
},K^{\varepsilon }\right)_{\varepsilon\in]0,1]}$
\end{lemma}
\begin{proof}
{\bf Step 1}: In view of step 1 of the proof of Lemma $\ref{C5L1}$
\begin{equation}
\sup_{\varepsilon }\sup_{n}\E\left(\sup_{0\leq s\leq t}\left|
Y_{s}^{\varepsilon ,n}\right| ^{2}+\int_{0}^{t}\left|
Z_{r}^{\varepsilon ,n}\right| ^{2}dr+\int_{0}^{t}\left|
Y_{r}^{\varepsilon ,n}\right| ^{2}dG_{r}^{\varepsilon }+\left|
K_{t}^{\varepsilon ,n}\right| ^{2}\right)<+\infty. \label{5.13}
\end{equation}
{\bf Step 2}: There exists an $\mathcal{F}_{s}$-adapted process $(\overline{Y}^{\varepsilon}_s)_{0\leq s\leq t}$ such that
\begin{eqnarray*}
\lim_{n\rightarrow \infty}\E\int_{0}^{t}|Y^{\varepsilon,n}_s-\overline{Y}^{\varepsilon}_s|^{2}ds=0.
\end{eqnarray*}
Notice that, for all $n\geq 1$ and $(s,x, y, z)\in[0,t]\times\R^{d}\times\R\times\R^{d}$,
\begin{eqnarray*}
f_n(s,x,y,z)\leq f_{n+1}(s,x,y,z).
\end{eqnarray*}
Therefore, by the comparison theorem for generalized BSDE \cite{PZ}, we have $Y^{\varepsilon,n}_s\leq Y^{\varepsilon,n+1}_s,\, a.s.$ for all $0\leq s\leq t$.
Hence, for all $0\leq s\leq t$
\begin{eqnarray*}
Y^{\varepsilon,n}_s\nearrow \overline{Y}^{\varepsilon}_s,\;\;\; a.s.
\end{eqnarray*}
In view of inequality $(\ref{5.13})$ and Fatou's lemma,
\begin{eqnarray*}
\E\left(\sup_{0\leq s\leq t}\left| \overline{Y}_{s}^{\varepsilon }\right|
^{2}\right)<C.
\end{eqnarray*}
It then follows by the Lebesgue convergence theorem that
\[
\lim_{n}\E\int_{0}^{t}|Y_{r}^ {\varepsilon
,n}-\overline{Y} _{r}^{\varepsilon }|^{2}dr=0.
\]
{\bf Step 3}: $\lim_{n\rightarrow}\E(\sup_{0\leq s\leq t}\left| \left(Y_{s}^{\varepsilon
,n}-h\left(s,X_{s}^{\varepsilon }\right)\right)^{-}\right| ^{2})=0$ \; a.s.

Let  $\left\{(\widetilde{Y}_{s}^{\varepsilon,n}\widetilde{Z}
_{s}^{\varepsilon,n}),%
\text{ }0\leq s\leq t\right\} $ be the unique solution of BSDE
\begin{eqnarray*}
\widetilde{Y}_{s}^{\varepsilon ,n} &=&l(X_{t}^{\varepsilon
})+\int_{s}^{t}f(X_{r}^{\varepsilon },Y_{r}^{\varepsilon
,n})dr+\int_{s}^{t}g(X_{r}^{\varepsilon },Y_{r}^{\varepsilon
,n})dG_{r}^{\varepsilon } \\
&&+n\int_{s}^{t}(h\left(r,X_{r}^{\varepsilon }\right)-\widetilde{Y}_{r}^{\varepsilon
,n})dr -\int_{s}^{t}\widetilde{Z}
_{r}^{\varepsilon ,n}dM^{X^{\varepsilon}}_{r}
\end{eqnarray*}
By using the comparison theorem, it follow for all $0\leq s\leq t$ and $n\geq 1$ that $\widetilde{Y}_{s}^{\varepsilon
,n}\leq Y_{s}^{\varepsilon ,n}$, a.s.

Now let $\nu$ be an $\mathcal{F}_{t}$-stopping time such that $0\leq \nu \leq t$. Then, applying Itô's formula to $\widetilde{Y}_{s}^{\varepsilon,n} e^{n(t-\nu)}$, we have
\begin{eqnarray*}
\widetilde{Y}_{\nu }^{\varepsilon,n} &=&\E^{\mathcal{F}_{\nu
}}[e^{-n\left(t-\nu \right)
}l(X_{t}^{\varepsilon })+\int_{\nu }^{t}e^{-n\left( s-\nu \right)
}f(X_{r}^{\varepsilon },Y_{r}^{\varepsilon ,n})dr \\
&&+\int_{\nu }^{t}e^{-n\left(
r-\nu \right)}g(X_{r}^{\varepsilon
},Y_{r}^{\varepsilon ,n})dG_{r}^{\varepsilon
} +n\int_{\nu }^{t}e^{-n\left(
r-\nu \right)}h\left(r,X_{r}^{\varepsilon }\right)dr].
\end{eqnarray*}
It is easily seen that
\begin{eqnarray*}
 e^{-n\left(t-\nu\right) }l(X_{t}^{\varepsilon })+n\int_{\nu }^{t}e^{-n\left(r-\nu \right)
}h\left(r,X_{r}^{\varepsilon }\right)dr\rightarrow l(X_{t}^{\varepsilon }) \mathbf{1}_{\left\{ \nu =t\right\}
}+h\left(t,X_{\nu}^{\varepsilon }\right)\mathbf{1}_{\left\{ \nu
<t\right\}}
\end{eqnarray*}
 a.s. and in $L^{2}\left( \Omega \right)$, and the conditional expectation converges also in $L^{2}\left( \Omega \right)$. Moreover,
\begin{eqnarray*}
\left|\int_{\nu }^{t}e^{-n\left( r-\nu \right)}f(X_{r}^{\varepsilon},Y_{r}^{\varepsilon ,n})dr\right|&\leq&\frac{1}{\sqrt{n}}\left(\int_{0}^{t}|f(X_{r}^{\varepsilon},Y_{r}^{\varepsilon ,n})|^{2}dr\right)^{1/2}\\
\left|\int_{\nu }^{t}e^{-n\left( r-\nu \right)}g(X_{r}^{\varepsilon},Y_{r}^{\varepsilon ,n})dG^{\varepsilon}r\right|&\leq&\frac{1}{\sqrt{n}}\left(\int_{0}^{t}|g(X_{r}^{\varepsilon},Y_{r}^{\varepsilon ,n})|^{2}dG^{\varepsilon}_{r}\right)^{1/2},
\end{eqnarray*}
and
\begin{eqnarray*}
\E^{\mathcal{F}_{\nu }}\left\{ \int_{\nu }^{t}e^{-n\left( r-\nu \right)
}f(X_{r}^{\varepsilon
},Y_{r}^{\varepsilon ,n})dr+\int_{\nu }^{t}e^{-n\left( r-\nu \right)}g(X_{r}^{\varepsilon
},Y_{r}^{\varepsilon ,n})dG_{r}^{\varepsilon }\right\}
\rightarrow 0
\end{eqnarray*}
in $L^{2}\left( \Omega \right)$ as $n\rightarrow \infty$.
Consequently,
\begin{eqnarray*}
\widetilde{Y}_{\nu }^{\varepsilon ,n} &\longrightarrow
&l(X_{t}^{\varepsilon })\mathbf{1}_{\left\{ \nu =t\right\}
}+h\left(\nu,X_{\nu}^{\varepsilon } \right)\mathbf{1}
_{\left\{ \nu <t\right\} } \;\mbox{in}\; L^{2}\left( \Omega \right)\mbox{as}\;
n \rightarrow\infty.
\end{eqnarray*}
Therefore, $\overline{Y}_{\nu}\geq h(\nu,X_{\nu})$ a.s.\newline
From this and the section theorem (see Dellacherie et Meyer \cite{DM}
page $220$), the previous inequality hold for all $s\in[0,t]$ a.s. Then
\begin{eqnarray*}
(Y^{n}_s-h(s,X_s))^{-}\searrow 0, \, 0\leq s\leq t,\; \mbox{a.s.},
\end{eqnarray*}
and from Dini's theorem the convergence is uniform in $t$. Since $(Y^{n}_s-h(s,X_s))^{-}\leq (h(s,X_s)-Y^{0}_s)^{+}\leq |h(s,X_s)|+|\overline{Y}^{0}_s|$, the result follows by dominated convergence.

{\bf Step 4}:\;$\lim_{n\rightarrow}\E\left(\sup_{0\leq s\leq t}|Y^{\varepsilon,n}_s-\overline{Y}_s|^{2}\right)=0$ and there exists $\mathcal{F}_{t}$-adapted processes $(\overline{Z}_s)_{0\leq s\leq t}$ and $(\overline{K}_s)_{0\leq s\leq t}$ such that
\begin{eqnarray*}
\lim_{n\rightarrow}\E\left(\int^{T}_{0}|Z^{\varepsilon,n}-\overline{Z}|^{2}+\sup_{0\leq s\leq t}|K^{\varepsilon,n}_s-\overline{K}_s|^{2}\right)=0.
\end{eqnarray*}

Indeed, by Itô's formula and taking expectation in both sides, we have
for all $n\geq m \geq 1$ and $0\leq s\leq t$,
\begin{eqnarray}
&&\left| Y_{s}^{\varepsilon ,n}-Y_{s}^{\varepsilon ,m}\right|
^{2}+\int_{s}^{t}\left| Z_{r}^{\varepsilon ,n}-Z_{r}^{\varepsilon
,m}\right| ^{2}dr+2\left| \beta \right| \int_{s}^{t}\left|
Y_{r}^{\varepsilon
,n}-Y_{r}^{\varepsilon ,m}\right| ^{2}dG_{r}^{\varepsilon }  \nonumber \\
&\leq &2\mu \int_{s}^{t}\left| Y_{r}^{\varepsilon
,n}-Y_{r}^{\varepsilon ,m}\right|^{2}dr
+2\int_{s}^{t}(Y_{r}^{\varepsilon ,n}-h\left(t,X_{r}^{\varepsilon }
\right))^{-}dK_{r}^{\varepsilon ,n}  \nonumber \\
&&+2\int_{s}^{t}(Y_{r}^{\varepsilon ,m}-h\left(r,X_{r}^{\varepsilon }\right))^{-}dK_{r}^{\varepsilon
,m}-\int_{s}^{t}\left(Y_{r}^{\varepsilon ,n}-Y_{r}^{\varepsilon
,m}\right)(Z_{r}^{\varepsilon ,n}-Z_{r}^{\varepsilon ,m})dM^{X^{\varepsilon}}_{r}.
\label{5.19}
\end{eqnarray}
Therefore
\begin{eqnarray*}
\E\int_{s}^{t}\left| Z_{r}^{\varepsilon ,n}-Z_{r}^{\varepsilon
,m}\right| ^{2}dr &\leq &C\E\int_{s}^{t}\left| Y_{r}^{\varepsilon
,n}-Y_{r}^{\varepsilon ,m}\right| ^{2}dr+2\E\int_{s}^{t}(Y_{r}^{%
\varepsilon ,n}-h\left(r,X_{r}^{\varepsilon }\right))^{-}dK_{r}^{\varepsilon ,n} \nonumber
\\
&&+2\E\int_{s}^{t}(Y_{r}^{\varepsilon
,m}-h\left(r,X_{r}^{\varepsilon }\right))^{-}dK_{r}^{\varepsilon
,m}.
\end{eqnarray*}
Furthermore, Hölder's inequality provide
\begin{eqnarray}
\E\int_{s}^{t}(Y_{r}^{\varepsilon ,n}-h\left(r,X_{r}^{\varepsilon }\right))^{-}dK_{r}^{\varepsilon
,n}\leq \left[\E\left(\sup_{0\leq s\leq t}|(Y^{\varepsilon,n}_s-h(s,X_s))^{-}|^{2}\right)\right]^{1/2}[\E|K_t^{m}|^{2}]^{1/2}
\label{Inq1}
\end{eqnarray}
and
\begin{eqnarray}
\E\int_{s}^{t}(Y_{r}^{\varepsilon ,p}-h\left(r,X_{r}^{\varepsilon }\right))^{-}dK_{r}^{\varepsilon
,n}\leq \left[\E\left(\sup_{0\leq s\leq t}|(Y^{\varepsilon,m}_s-h(s,X_s))^{-}|^{2}\right)\right]^{1/2}[\E|K_t^{n}|^{2}]^{1/2}.
\label{Inq2}
\end{eqnarray}
By using the result of Step 1 and Step 2, the right hands sides of the above inequalities $(\ref{Inq1})$ and $(\ref{Inq2})$ go to zero as $m,n \rightarrow \infty$. It follows that
\begin{eqnarray*}
\E\int_{0}^{t}\left| Z_{r}^{\varepsilon ,n}-Z_{r}^{\varepsilon
,m}\right| ^{2}dr \rightarrow  0
\end{eqnarray*}
as $n, m \rightarrow \infty$, and there exists a process $\overline{Z}^{\varepsilon}\in H^{2,d}(0,t)$ such that
\begin{eqnarray*}
\E\int_{0}^{t}\left| Z_{r}^{\varepsilon ,n}-\overline{Z}_{r}^{\varepsilon}\right| ^{2}dr \rightarrow  0,\; \mbox{as}\; n\rightarrow \infty.
\end{eqnarray*}
Then, by the Burkholder-Davis-Gundy inequality, we get
\begin{eqnarray*}
\E(\sup_{0\leq s\leq t}\left| Y_{s}^{\varepsilon,n}-Y_{s}^{\varepsilon ,m}\right| ^{2}&\leq &4\E\int_{0}^{t}(Y^{n}-S_s)^{-}dK^{m}_s+4\E\int_{0}^{t}(Y^{m}-S_s)^{-}dK^{n}_s\\
&&+C\E\int_{0}^{t}\left| Z_{r}^{\varepsilon ,n}-Z_{r}^{\varepsilon
,m}\right|ds.
\end{eqnarray*}
It follows that $\E(\sup_{0\leq s\leq t}\left| Y_{s}^{\varepsilon,n}-Y_{s}^{\varepsilon ,m}\right| ^{2}\rightarrow\,\mbox{as}\, n,m \rightarrow \infty$. Therefore,
\begin{eqnarray*}
\E\left\{ \sup_{0\leq s\leq t}\left| Y_{s}^{\varepsilon
,n}-\overline{Y}_{s}^{\varepsilon}\right| ^{2}\right\}\rightarrow 0\,\mbox{as}\, n\rightarrow \infty.
\end{eqnarray*}
We deduce that
\begin{eqnarray*}
\E\left(\sup_{0\leq s\leq t}\left| K_{s}^{\varepsilon ,n}-K_{s}^{\varepsilon ,m}\right| ^{2}\right)
\rightarrow 0\, \mbox{as}\, n,m \rightarrow \infty.
\end{eqnarray*}
Hence there exists a process
$(\overline{K}_{s}^{\varepsilon })_{0\leq s\leq t}\in A^{2}(0,t)$ such that
\begin{eqnarray*}
\E\left(\sup_{0\leq s\leq t}\left| K_{s}^{\varepsilon ,n}-\overline{K}_{s}^{\varepsilon}\right| ^{2}\right)
\rightarrow 0\, \mbox{as}\, n\rightarrow \infty.
\end{eqnarray*}

{\bf Step 5}:\; The limiting process $(\overline{Y}^{\varepsilon},\overline{Z}^{\varepsilon},\overline{K}^{\varepsilon})$ is the solution of the reflected generalized BSDE $(\xi,f,g,h)$.\\
Obviously the process $(\overline{Y}^{\varepsilon}_s,\overline{Z}
^{\varepsilon}_s,\overline{K}^{\varepsilon}_s)_{0\leq s\leq t}$ satisfies
\begin{eqnarray*}
\overline{Y}^{\varepsilon}_s=\xi + \int^{t}_{s}f(r,X^{\varepsilon}_r,\overline{Y}^{\varepsilon}_r,\overline{Z}^{\varepsilon}_r)ds + \int^{t}_{s}
g(r,X^{\varepsilon}_r,\overline{Y}^{\varepsilon}_r)dG^{\varepsilon}_r + \overline{K}^{\varepsilon}_t - \overline{K}^{\varepsilon}_s - \int^{t}_{s}
\overline{Z}^{\varepsilon}_r dM^{X^{\varepsilon}}_r,\,\, 0\leq s\leq t.
\end{eqnarray*}
Since $(Y^{\varepsilon,n},K
^{\varepsilon,n})$ tends to $\left(
\overline{Y}^{\varepsilon },\overline{K}^{\varepsilon}\right) $ in probability uniformly in $t$ ,
the measure $dK^{\varepsilon,n}$ converges to $d\overline{K}^{\varepsilon}$ weakly in probability, so that
$\int_{0}^{t}(Y_{r}^{\varepsilon ,n}-h\left(r,X_{r}^{\varepsilon }\right))dK_{r}^{\varepsilon,n} \rightarrow \int_{0}^{t}(\overline{Y}_{r}^{\varepsilon }-h\left(r,X_{r}^{\varepsilon }\right))d\overline{K}^{\varepsilon }_{r}$ as $n \rightarrow\infty$. Obviously, $\int_{0}^{t}(\overline{Y}_{r}^{\varepsilon }-h\left(r,X_{r}^{\varepsilon }\right))d\overline{K}^{\varepsilon }_{r}\geq 0$, while, on the other hand, for all $n\geq 0$, $\int_{0}^{t}(Y_{r}^{\varepsilon ,n}-h\left(r,X_{r}^{\varepsilon }\right))dK_{r}^{\varepsilon,n}\leq 0$.\\
Hence,
\begin{eqnarray*}
\int_{0}^{t}\left( \overline{Y}_{s}^{\varepsilon
}-h\left(s,X_{s}^{\varepsilon }\right)\right)
d\overline{K}_{s}^{\varepsilon }=0.\;\;\;\; \mbox{a.s}.
\end{eqnarray*}
Consequently, $(\overline{Y}_{s}^{\varepsilon },\overline{Z}_{s}^{\varepsilon },\overline{K}_{s}^{\varepsilon })_{0\leq s\leq t}$ is the solution of the reflected generalized BSDE
$(\xi,f,g,h)$ and uniqueness of this equation allows that $\overline{Y}^{\varepsilon }=Y^{\varepsilon },$
$\overline{Z} ^{\varepsilon }=Z^{\varepsilon }$ and
$\overline{K}^{\varepsilon }=K^{\varepsilon }$ which ends the proof.
\end{proof}
\newline\newline\newline
{\bf Proof of Theorem $\ref{C5T1}$}\\
Put Theorem $2.1$, Lemma $\ref{C5L1},$ Lemma $\ref{C5L3}$, Lemma
$\ref{C5L2}$ and the Proposition $\ref{C5P1}$ together prove that
$\left(Y^{\varepsilon},M^{\varepsilon },K^{\varepsilon }\right) $ converges in law to
$\left(Y,M,K\right) $ in $D([0,T],\R)\times D([0,T],\R)\times C([0,T],\R_+)$, where the first factor is equipped with the $S$-topology, and the two last
factors with the topology of uniform convergence. In particular
 \[
Y_{0}^{\varepsilon }=\E\left\{ l\left( X_{t}^{\varepsilon }\right)
+\int_{0}^{t}f\left(X_{s}^{\varepsilon },Y_{s}^{\varepsilon
}\right) ds+\int_{0}^{t}g\left(X_{s}^{\varepsilon
},Y_{s}^{\varepsilon }\right) dG_{s}^{\varepsilon
}+K_{t}^{\varepsilon}\right\} .
\]
converge to
 \[
Y_{0}=\E\left\{ l\left( X_{t}\right)
+\int_{0}^{t}f\left(X_{s},Y_{s}\right)ds+\int_{0}^{t}g\left(X_{s},Y_{s}\right)dG_{s}+K_{t}\right\} .
\]
in $\R$.
Indeed, since $Y_{0}^{\varepsilon }$ is deterministic,
$Y_{0}^{\varepsilon }=\E(B_{\varepsilon})$ where
\[
B_{\varepsilon }=l\left( X_{t}^{\varepsilon }\right)
+\int_{0}^{t}f\left(X_{s}^{\varepsilon },Y_{s}^{\varepsilon
}\right) ds+\int_{0}^{t}g\left(X_{s}^{\varepsilon
},Y_{s}^{\varepsilon }\right) dG_{s}^{\varepsilon
}+K_{t}^{\varepsilon }.
\]
In view of  assumptions $({\bf H1})$-$({\bf H4})$, we have,
\begin{eqnarray*}
\sup_{\varepsilon}\E| B_{\varepsilon}|^{2}\leq C.
\end{eqnarray*}
Since $B_{\varepsilon}$ converge in law as $\varepsilon$ goes to $0$, toward
\[B=l\left( X_{t}\right) +\int_{0}^{t}f\left(X_{s},Y_{s}\right)
ds+\int_{0}^{t}g\left(X_{s},Y_{s}\right) dG_{s}+K_{t},
\]
the uniformly integrability of $B_{\varepsilon }$ implies that
\[
\lim_{\varepsilon\rightarrow 0}\E\left( B_{\varepsilon }\right)=\E\left(\lim_{\varepsilon\rightarrow 0}B^{\varepsilon}\right).
\]
This mean that $Y^{\varepsilon}$ converges to $Y_0=\E(B)$.\\
$\blacksquare $

\section {Homogenization of the reflected viscosity solutions
of scalar PDE with non linear Neumann boundary condition}
\setcounter{theorem}{0} \setcounter{equation}{0}

The result of the above section permits to us to deduce weak convergence of a sequence $(Y^{\varepsilon})$ of the solutions of reflected generalized BSDEs from weak convergence of the sequence $(X^{\varepsilon})$. In a sense, we deduce, from the probabilistic proof of the convergence of linear
PDE with linear boundary condition due to Tanaka \cite{T}, a probabilistic proof of convergence of reflected semi-linear PDEs with non linear Neumann boundary condition.

For each $x\in\bar{\Theta}$, let $\{(X^{\varepsilon,x}_s,G^{\varepsilon,x}_s),\ s\geq 0\}$ denote the solution of the SDE $(\ref{a4})$. Let $l, f,g, h$ be as in the previous section. For each $(t,x)\in\R_+\times\bar{\Theta}$, let
\[
u^{\varepsilon}(t,x)=Y^{\varepsilon,x}_0,
\]
where $Y^{\varepsilon}$ denote the solution of the reflected generalized BSDE $(\ref{a6})$ considered in the previous section (which of course depends on the starting point $x$ of $X^{\varepsilon}$, and the final time $t$ for the reflected generalized BSDE). The continuity of $u^{\varepsilon}$ follows by the standard arguments (see \cite{Ral}) from the continuity of the mapping $x\mapsto X^{\varepsilon,x}$ in probability. One can show as in Ren and Xia \cite{Ral} that $u^{\varepsilon}$ is the viscosity solution of the semi-linear parabolic PDEs
\begin{eqnarray}
\left\{
\begin{array}{l}
 \min \left(u^{\varepsilon }\left( t,x\right) - h\left( t,x\right),\right. \\\\
\left.-\frac{\partial u^{\varepsilon }}{\partial t}\left( t,x\right)
- L^{\varepsilon }_{x}u^{\varepsilon }\left( t,x\right)
-f\left(x,u^{\varepsilon }\left( t,x\right) \right) \right)=0,\,\, if\,\, \left( t,x\right)\in \R_{+} \times \Theta , \\\\
\Gamma ^{\varepsilon }u^{\varepsilon }\left( t,x\right)
+g\left(x,u^{\varepsilon }\left( t,x\right) \right) =0,\,\, if \,\, \left( t,x\right)
\in \R_{+} \times \partial \Theta , \\\\
u^{\varepsilon }\left( 0,x\right) =l\left( x\right) ,\text{ }x\in \overline{%
\Theta },
\end{array}\right.  \label{5.23}
\end{eqnarray}
where $L_{x}^{\varepsilon}$ and $\Gamma^{\varepsilon}_x$ are the operators defined by $(\ref{a2.})$ and $(\ref{a2})$ respectively.

Let $u$ be the solution of the homogenized PDE
\begin{eqnarray}
\left\{
\begin{array}{l}
 \min\left(u\left( t,x\right) - h\left( t,x\right),\right. \\\\
\left.-\frac{\partial u}{\partial t}\left( t,x\right) - L^{0}_{x}u \left(
t,x\right) -f\left(x,u\left( t,x\right) \right) \right)=0,\,\, if \,\,\left(
t,x\right) \in \R_{+} \times \Theta , \\\\
\Gamma u\left( t,x\right) +g\left(x,u\left( t,x\right) \right) =0,\,\, if\,\,
\left( t,x\right) \in \R_{+} \times \partial \Theta,
\\\\
u\left( 0,x\right) =l\left( x\right) ,\text{ }x\in
\overline{\Theta },
\end{array}
\right.  \label{5.24}
\end{eqnarray}
where $L^{0}_{x}$ and $\Gamma^{0}_{x}$ are defined in $(\ref{op})$. Recall again the previous argument, $u(t,x)=Y_{0}^{x}$ is the viscosity solution of $(\ref{5.24})$, where $Y^{x}$ denote the solution of the reflected generalized BSDE $(\ref{Eq1})$, which depend on the starting point of $X^{x}$.
\begin{theorem}
\label{C4T2}  Under the assumptions $(\bf H1)$-$(\bf H4)$, for each $\left( t,x\right)\in \R_{+} \times \Theta$, Then  $
u^{\varepsilon }\left( t,x\right)$ converges to $u\left(
t,x\right) $ as $\varepsilon$ goes to 0.
\end{theorem}

\begin{proof}
Since $u^{\varepsilon}(t,x)=Y^{\varepsilon,x}_{0}$ and $u(t,x)=Y^{x}_{0}$, the result follows from Theorem $3.1$
\end{proof}

\end{document}